\documentclass{article}
\usepackage{amsmath}
\usepackage{amssymb}
\usepackage{amscd}
\usepackage[matrix]{xypic}
\usepackage{cite}

\newtheorem{theorem}{Theorem}

\newtheorem{corollary}[theorem]{Corollary}

\newtheorem{definition}[theorem]{Definition}

\newtheorem{proposition}[theorem]{Proposition}

\newcommand\p{\mathcal{P}}

\begin{document}

\title{On the nonKoszulity of $(2p+1)$-ary partially associative Operads}
\author{Elisabeth REMM \\
elisabeth.remm@uha.fr\\
Adress : LMIA.\ UHA\\
4 rue des Fr\`{e}res Lumi\`{e}re\\
F. 68093 Mulhouse Cedex}
\maketitle

\begin{abstract}
We want to present here the part of the work in common with Martin Markl \cite{M.R} which concerns quadratic operads  
for $n$-ary algebras and their dual for $n$ odd. We will focus on the ternary case (i.e $n=3$). 
The aim is to underline the problem of computing the dual operad and the fact that this last is in general defined
in the graded differential operad framework. 

We prove that the operad associated to $(2p+1)$-ary partially associative algebra is not
Koszul. Recall that, in the even case, this operad is Koszul.
\end{abstract}

\section{\protect\bigskip Introduction}

We are interested in the operads associated to $n$-ary algebras
whose multiplication $\mu $ satisfies the following relation%
\[
\underset{i=1}{\overset{n}{\sum }}(-1)^{(i-1)(n-1)}\mu (X_{1},\cdots ,X_{i-1},\mu
(X_{i},\cdots ,X_{i+n-1}),X_{i+n},\cdots ,X_{2n-1})=0.
\]
Such a multiplication is called $n$-ary partially associative.

These operads have already been studied several times for instance in the
thesis of Gnedbaye \cite{Gne} and more recently by M. Markl - E. Remm (in term of
operads) \cite{M.R}, N. Goze - E. Remm (in term of algebras) \cite{G.R} and also H. Ataguema - A. Makhlouf \cite{At.Mk} and 
E. Hoffbeck \cite{H}.

When computing the free algebras associated to $(2p+1)$-ary partially
associative algebras, using different arguments in  \cite{M.R} and  \cite{G.R}, 
we saw that the odd  and even cases
behave in a completely different way.
This approach is a little bit different to  \cite{Gne}, \cite{At.Mk} and \cite{G.W} where
odd and even cases have been studied in the same way and then contain some
  misunderstandings in the odd case. 
In fact, contrary to what we find in the
previous papers, the operads associated to $(2p+1)$-ary partially
associative algebras are non Koszul so there is no operadic cohomology
associated to these algebras (we call it operadic because it uses the dual operad to define a cohomology). 
It is this result that we are going to expose
in this work.

Let us remark that in the first versions of \cite{H}, the author thought that the degrees didn't matter and that
 the degrees of the operations
appear only in the calculation of the dual and even
there, if the operations are concentrated in a fixed arity, it wouldn't
change the relations in the dual. Then he concluded that the operad of partially associative $n$-ary algebras was Koszul for
any $n$.  The author has changed 
his point of view after that we indicated him the problem. 
Let us also note that in his habilitation thesis and contrary to his previous works Gnedbaye 
precises the necessity of distinguishing the even and odd cases
 assigning the problem on the definition of the generating series 
what is is insufficient to find the odd case as we are going to see here.

To study Koszulity we need to define correctely the dual operad of a $(2p+1)$%
-ary partially associative operad which implies to understand correctely the definition of
Ginzburg-Kapranov \cite{G.K} developed in the case of binary operations in order to
extend it to $n$-ary operations. When we compute the dual of a graded or nongraded 
operad, we get graded objects which take some supension in account. 
In particular if we consider  an $n$-ary multiplication of degree $0$ and its associated operad $\mathcal{P}$, 
an algebra on its dual operad $\mathcal{P}^!$ corresponds 
to a $n$-ary multiplication of degree $m\equiv n \ modulo \ 2$. 
In our case, as $n$ is even, the multiplication of an algebra on the dual operad 
has to be considered of degree $1$. 
If we forgot this degree, all multiplications related to the dual operad are considered of degree $0$. We get something
that we could call the nongraded "dual" operad which is the operad of $(2p+1)$-ary totally associative
algebras with operation in degree $0$. But this operad is not the dual operad in the Ginzburg-Kapranov's sense
when $n$ is odd.

As it was not possible to find an operadic  cohomology for $(2p+1)$-ary partially associative algebras, 
 we have developped in \cite{G.R} a cohomology for $(2p+1)$-ary partially associative algebras  
restricting the space of cochains that we developp for the ternary case. We give explicitely the free algebra.  
In \cite{M.R}
we have a different approach. We consider a \textbf{graded}
version of the $(2p+1)$-ary partially associative algebras with operation in degree $1$. In this case the dual operad
is the operad of $(2p+1)$-ary totally associative algebras with operation in degree $0$ and both are Kozul.
Then we obtain an operadic cohomology.

In the following we consider $\mathbb{K}$ a field of characteristic $0$ and the operads that we consider
are generally $\mathbb{K}$-linear operads. All definitions and concepts used refers to \cite{G.K} and \cite{M.S.S}.

\section{The operad $3\mathcal{A}ss$}

This section deals with the operad of $3$-ary partially associative algebras that is algebras defined by a multiplication
$$\mu :A^{\otimes 3}\rightarrow A$$
satisfying the relations
$$\mu \circ (\mu \otimes I_2)+\mu \circ (I_1 \otimes \mu \otimes I_1)+
\mu \circ (I_2 \otimes \mu)=0$$
where $I_j:A^{\otimes j}\rightarrow A^{\otimes j}$ is the identity map. We have a classical example of such an algebra
when we consider the Hochschild cohomology of an associative algebra. In fact if 
$\mathcal{C}^k(V,V)$ denotes the space of $k$-cochains of the Hochschild cohomology of the 
associative algebra $V$, the Gerstenhaber product $\circ_{n,m}$ is a linear map
$$\circ_{n,m}: \mathcal{C}^n(V,V) \times \mathcal{C}^m(V,V) \rightarrow \mathcal{C}^{n+m-1}(V,V)$$
given by
$$
\begin{array}{l}
f\circ _{n,m}g(X_{1}\otimes \cdots \otimes
X_{n+m-1})=\\
\\
\displaystyle\sum_{i=1}^{n}(-1)^{(i-1)(m-1)}f(X_{1}\otimes \cdots \otimes
g(X_{i}\otimes \cdots \otimes X_{i+m-1})\otimes \cdots \otimes X_{n+m-1})
\end{array}
$$
if $f\in C^{n}(V)$ and $g\in C^{m}(V).$ 
A $3$-ary partially associative product is a $3$-cochain $\mu$ satisfying 
$\mu \circ_{3,3} \mu =0.$

The notion of quadratic operad is clearly defined in [G.K]. We refer to this paper. An operad $\mathcal{P}$
is a collection  $\left\{ \mathcal{P}(n) , n\geq 1 \right\}$ of $\mathbb{K}$-vector spaces 
satisfying the following properties: each $\mathcal{P}(n)$ is a $\Sigma _n$-module where 
$\Sigma _n$ is the symmetric group of degree $n,$ there is an element $1 \in \mathcal{P}(1)$ called the unit 
and linear maps
$$\circ_i:  \mathcal{P}(n) \times  \mathcal{P}(m) \rightarrow  \mathcal{P}(n+m-1)$$
called comp-i operations satisfying associative conditions:
if $\lambda \in  \mathcal{P}(l), \mu \in  \mathcal{P}(m), \nu  \in  \mathcal{P}(n)$ then
$$(\lambda \circ_i \mu )\circ_{j} \nu =
\left\{ 
\begin{array}{l}
(\lambda \circ_j \nu )\circ_{i+n-1} \mu  \ \mbox{\rm  if} \ 1\leq j\leq i-1 \\
\lambda \circ_i (\mu \circ_{j-i+1} \nu)  \ \mbox{\rm  if} \ i \leq j \leq m+1-1 \\
(\lambda \circ_{j-m+1} \nu )\circ_{i} \mu \ \mbox{\rm  if} \ i+m \leq j \\
\end{array}
\right.
$$
which are compatible with the action of the symmetric group.

Recall that 
a $\mathcal{P}$-algebra is a $\mathbb{K}$ vector space $V$ equipped with a morphism of operad 
 $f: \mathcal{P} \rightarrow \mathcal{E}_V$ where $\mathcal{E}_V$ is the operad of endomorphisms of $V$.
 Given a structure of $\mathcal{P}$-algebra on $V$ is the same as giving a collection of linear maps
$$f_n:\mathcal{P}(n)\otimes V^{\otimes n}\rightarrow V$$ satisfying natural 
associativity, equivariance and unit conditions.

If $E$ is  a right-$\mathbb{K}[\Sigma _2]$-module, we can define an operad, denoted by 
$\mathcal{F}(E)$ and called the free operad generated by $E$ which is solution of the 
following universal problem: for any operad $\mathcal{Q}=\{ \mathcal{Q}(n) \}$
and any $\mathbb{K}[\Sigma _2]$-linear morphism $f:E \rightarrow \mathcal{Q}(2)$, there exits
a unique operad morphism $\hat{f}: \mathcal{F}(E)\rightarrow \mathcal{Q}$ which coincide
with $f$ on $E=\mathcal{F}(E)(2).$ We have for example  
$\mathcal{F}(E)(3)=(E \otimes E) \otimes_{\Sigma _2} \mathbb{K}[\Sigma _3].$
If $R$ is a $\mathbb{K}[\Sigma _3]$-submodule of $\mathcal{F}(E)(3)$, it generates
an ideal $(\mathcal{R})$ of  $\mathcal{F}(E).$ The quadratic operad generated by $E$
with relations $R$ is the operad $\mathcal{P}(E,R)=\{ \mathcal{P}(E,R)(n),n\geq 1  \}$ with
$$\mathcal{P}(E,R)(n)=\mathcal{F}(E)(n)/\mathcal{R}(n).$$
This notion of quadratic operad is related to binary algebras.
In [Gn] this notion is adapted to $n$-ary algebras. 
In this case we consider a generating multiplication $\mu$  which is a $n$-ary multiplication that is 
$E=<\mu>$ a $ \mathbb{K}[\Sigma _n]$-module. We define the free operad generated by $E$ in the same sense that in the 
binary case. Here we get that $\mathcal{F}(E)(n)=0$ if $n \neq p(n-1)+1$ and   
$\mathcal{F}(E)(p(n-1)+1)$ consists as a vector space to 
"parenthesized products" of $p(n-1)+1$ variables indexed by $\left\{1, \cdots ,n \right\}$. For instance
a basis of $\mathcal{F}(E)(n)$ is generated as $\mathbb{K}[\sigma_n]$-module by $(x_1 \cdots x_n)$ and a basis 
of $\mathcal{F}(E)(2(n-1)+1)$ is given by 
$((x_1 \cdots x_n) x_{n+1} \cdots x_{2n-1}),(x_1 (x_2 \cdots x_{n+1}) x_{n+2} \cdots x_{2n-1}),\cdots ,$ 

\noindent $(x_1 \cdots x_{n-1}( x_{n} \cdots x_{2n-1})$ and all their permutations. The relations that we consider will
be quadratic in the sense that we compose two $n$-ary multiplications so
 $R$ is a $\mathbb{K}[\Sigma _{2n-1}]$-submodule of $\mathcal{F}(E)(2n-1)$. 
The 
operad $\mathcal{R}=(R)$ is the ideal generated by $R$ so in particular $\mathcal{R}(k)=0$ for $k=1$ and 
$k \neq p(n-1)+1$. 
In other words 
$\mathcal{R}(p(n-1)+1)$ consists in all relations in $\mathcal{F}(E)(p(n-1)+1)$ induced by the relations $R$.

We recall here the notion of quadratic operad
of ternary algebras.

\begin{definition}
Let $E$ be a $\mathbb{K}[\Sigma _3]$-module and $\mathcal{F}(E)$ the free operad over $E$. If 
$R$ is a  $\mathbb{K}[\Sigma _5]$-submodule of $\mathcal{F}(E)(5)$ and if $(R)$ is the ideal
generated by $R$, then the ternary quadratic operad $\mathcal{P}(E,R)$ generated by $E$ and $R$ 
is the quotient 
$$\mathcal{F}(E)/(R).$$
Note that $\mathcal{F}(E)(n)=0$ when $n$ is even.
If $E\simeq \mathbb{K}[\Sigma _3]$ and $R$ is the $\mathbb{K}[\Sigma _5]$-submodule of $\mathcal{F}(E)(5)$ 
generated by the vectors
$$ (x_1x_2x_3)x_4x_5+  x_1(x_2x_3x_4)x_5+  x_1x_2(x_3x_4x_5)$$
then the corresponding quadratic operad is 
the ternary quadratic operad  
$$3\mathcal{A}ss=\mathcal{F}(E)/(R).$$
\end{definition}   

We know that for any operad $\mathcal{P}$, the spaces $\mathcal{P}(n)$ are related to the
free $\mathcal{P}$-algebras. In [G.R] we have studied the free partially associative algebra of order $3$
$\mathcal{L}_{3\mathcal{A}ss}(V)=\oplus \mathcal{L}^{2p+1}(V)$ 
on a vector space $V$. We have computed the dimensions of its homogeneous components, found a basis
 and a 
systematic method to write this basis. In particular we have, if $dim V=1$, 
$$dim \mathcal{L}^3(V)=1, \ dim \mathcal{L}^5(V)=2, \ dim \mathcal{L}^7(V)=4$$
$$dim \mathcal{L}^9(V)=5, \ dim \mathcal{L}^{11}(V)=6, \ dim \mathcal{L}^{13}(V)=7. $$
We deduce that 
$$dim ((3\mathcal{A}ss)(3))=dim \mathbb{K}[\Sigma _3]=6, \ dim ((3\mathcal{A}ss)(5))=2\, dim \mathbb{K}[\Sigma _5]=240$$
and more generally
$$dim ((3\mathcal{A}ss)(2k+1))=(k+1)\, dim \mathbb{K}[\Sigma _{2k+1}].$$
The Poincar\'e serie of $(3\mathcal{A}ss)$ called also the generating function is written
$$g_{3\mathcal{A}ss}(x)=\sum_{n=1}^\infty dim((3\mathcal{A}ss)(n))\frac{x^n}{n!}$$
with the convention $dim (3\mathcal{A}ss)(1)=1$. Then
$$g_{3\mathcal{A}ss}(x)=x+x^3+2x^5+4x^7+5x^9+6x^{11}+7x^{13}+ \cdots$$

Recall that a quadratic operad $\p$ is Koszul if the homology of any free $\p$-algebra $F_{\p}(V)$ is trivial except in degree $0$.
If $\p$ is a Koszul operad, then its Poincar\'e serie $g_{\p}(x)$ satisfies 
$$g_{\p}(-g_{\p^!}(-x))=x$$
where $\p^!$ is the dual operad. In the following section we define the dual operad of $(3\mathcal{A}ss)$. 
The previous identity shows, if
$\p$ is of Koszul, that there is a formal serie $s(x)=\sum a_nx^n$ such that $g_{\p}(-s(-x))=x.$
Now consider a serie $s(x)=\sum a_nx^n$ and if we solve the equation $g_{3\mathcal{A}ss}(-s(-x))=x$ we find:
$$s(x)=x-x^3+x^5-19 x^{11}+O[x]^{12}.$$
Such a serie cannot be a Poincar\'e serie of a quadratic operad corresponding to a multiplication of degree $0$. 
In the next section we compute the graded dual operad and we will see that its generating function is a 
polynomial of degree $5$. 
This will permit to conclude to the non Koszulity of $(3\mathcal{A}ss)$. 

\medskip

\noindent {\bf Remark}. If we consider the operad $\mathcal{P}=3t\mathcal{A}ss$ associated to $3$-ary totally associative algebras
that is satisfying the relations
$$\mu \circ (\mu \otimes I_2)=\mu \circ (I_1 \otimes \mu \otimes I_1)=
\mu \circ (I_2 \otimes \mu)=0$$ we get the generating function
$$g_{3t\mathcal{A}ss}(x)=x+x^3+x^5+x^7+x^9+ \cdots + x^{2k+1}+ \cdots $$
If we suppose this operad to be Koszul we would get that the dual operad $\mathcal{P}^!$ has as generating serie
$g_{\mathcal{P}^!}$ satisfying $g_{\p}(-g_{\p^!}(-x))=x$. But if we compute $g_{\p^!}$ from these
equation we obtain $$g_{\p^!}(x)=\sum_{i \geq 1}a_n x^n$$   
but $\mid  a_n \mid $ do not correspond 
to the dimensions of the operad $3\mathcal{A}ss$.       

\section{The dual operad}

To compute the dual operad of an operad associated to $n$-ary algebras we need some differential {\bf graded} operad.
Recall that if $\mathcal{C}$ is a monoidal category, a $\Sigma $-module $A$ is represented by a sequence
of objects $\{A(n) \}_{n\geq 1}$ in $\mathcal{C}$ with a right-$\Sigma _n$-action on $A(n)$. 
Then an operad in a (strict) symmetric monoidal category $\mathcal{C}$ is a $\Sigma $-module 
$\mathcal{P}$ together with a family of structural morphisms 
satisfying some associativity, equivariance and unit conditions.

We  use the fact that $dgVect$ the category of differential graded vector spaces over the base field $\mathbb{K}$
(an object of $dgVect$ is a graded vector space together with a linear map $d$ (differential) of degree $1$ such that
$d^2=0;$ morphisms are linear maps preserving gradings and differentials)  
is a symmetric monoidal category and we can define an operad in this category $dg Vect$.

A differential {\bf graded} operad (or dg operad) is a differential graded $\Sigma $-module with 
an operad structure for which the operad structure maps are differential graded
 morphisms. 

A (non graded) operad can be seen as a differential graded operad considering trivial differentials and
nongraded objects are objects trivially graded.

Since a $dg$ operad $\mathcal{P}$ is itself a monoid in a symmetric monoidal category, the bar construction 
applies to $\mathcal{P}$, producing a $dg$ operad $\mathcal{B}(\mathcal{P})$. 
The linear dual of $\mathcal{B}(\mathcal{P})$ 
is a $dg$ operad denoted by $\mathcal{C}(\mathcal{P})$ and called the cobar complex of the operad $\mathcal{P}$. 
We also need the dual $dg$ operad $\mathcal{D}(\mathcal{P})$, which is just $\mathcal{C}(\mathcal{P})$
suitably regarded.
Quadratic operad are defined as having a presentation with generators and relations and  for them
the dual operad will also be quadratic. For quadratic operads there is a natural transformation of functor
from $\mathcal{D}(\mathcal{P})$ to $\mathcal{P}^!$ which is a quasi-isomorphism if $\mathcal{P}$ is Koszul.
 This concept is similar to the concept of quadratic dual and
Kozulness for associative algebras.

Any quadratic operad $\mathcal{P}$ associated to binary multiplications admits a dual 
operad which is also quadratic and denoted $\mathcal{P}^!.$ 
To define it we need to recall some definitions and notations.
Let  $E$ be a $\Sigma $-module. The dual 
$\Sigma $-module $E^\# =\{ E^\# (n) \}_{n\geq 1}$ 
is defined by $$E^\# (n)=Hom_{\mathbb{K}}(E(n),\mathbb{K})$$
and the $\Sigma _n$ representation on $E(n)$ determines a dual representation on $E^\# (n)$ by
$$(\lambda \cdot \sigma ,\mu ):=(\lambda ,\mu  \cdot \sigma ^{-1})$$
for $\mu  \in E(n), \lambda  \in E^\# (n)$ and $\sigma \in \Sigma_n.$ 
The Czech dual is the $\Sigma $-module $E^\vee  =\{ E^\vee (n) \}_{n\geq 1}$  with 
$$E^\vee (n)=E^\# (n) \otimes sgn_n.$$
Then the  quadratic dual operad is defined as a quotient of the free operad $\mathcal{F}(E^\vee )$
by relations othogonal (in some sense) to the relations defining the original quadratic operad $\mathcal{P}$. 
So if $\mathcal{Q}=\mathcal{P}(E,R)$ (i.e $E$ corresponds to the generators and $R$ to the relations) the dual operad
$\mathcal{Q}^!$ is defined by
$\mathcal{Q}^!=\mathcal{P}(E^\#  ,R^\bot  )$ where $R^\bot  \subseteq \mathcal{F}(E^\#  )(3)$
is the annihilator wih respect to some paring of the relations 
$R\subseteq \mathcal{F}(E)(3)$ defining $\mathcal{Q}$.
But notice that in the general definition of quadratic operad  contains a suspension. 
If we consider a quadratic operad generated by an operation 
$E=<\mu > \in \mathbb{K}[\Sigma _2]$ where $\mu$ is of degree $0$, the dual operad  
is still a quadratic operad generated by an operation  of degree $0$. 

Now if we consider $n$-ary algebras,  we have seen that we can still define the notion of quadratic operad
that is we consider a generating multiplication which is a $n$-ary multiplication $\mu$ that is 
$E=<\mu>\subset \mathbb{K}[\Sigma _n]$
and relations which are quadratic that is  $R$ is a $\mathbb{K}[\Sigma _{2n-1}]$-submodule. 
We can also still define a "scalar product" as in the case of a binary operation 
but now it is a map  $<\, , \, >: \mathcal{F}(E)(n) \otimes \mathcal{F}(E^\vee)(n) \rightarrow \mathbb{K}.$
Then we get $R^\perp $ but the dual operad $\mathcal{P}^!$ is 
$$\mathcal{P}^!(n)=\mathcal{F}(E^\vee )(n)/(R^\perp )(n)  $$  
when the generating operation is of degree even. But if the degree odd, the dual operad is a quadratic operad
with a {\it generating operation of degree odd} and other relations are induced by $R^\perp $.

This can also be seen as follows: if the operations are not of degree 0, then a nontrivial sign of the composition in
the free operad which is introduced.

Now let us come back to the determination to the dual operad of $3\mathcal{A}ss$. Recall that a totally associative 
$3$-ary algebra is given by a $3$-ary product $\mu$
satisfying
$$ \mu(\mu(x_1,x_2,x_3),x_4,x_5)=\mu(x_1,\mu(x_2,x_3,x_4),x_5)=\mu(x_1,x_2,\mu(x_3,x_4,x_5)).$$
We denote by $3p\mathcal{A}ss$ the corresponding quadratic operad. 

\begin{theorem}\cite{G.R}
The dual operad  of the $3\mathcal{A}ss$ operad is isomorphic to the operad of
totally associative algebras with operation of degree $1.$
\end{theorem}

Hence  as a vector space a basis of $3\mathcal{A}ss^{!}(n)$ is given by $( x_{i_1} x_{i_2} \cdots  x_{i_n})$,
a basis of   
$3\mathcal{A}ss^{!}(2n-1)$ is given by $(x_{i_1} \cdots x_{i_{n-1}}( x_{i_n} \cdots x_{i_{2n-1}}))$. We have also 
$$3\mathcal{A}ss(n)^{!}(k(n-1)+1)=\left\{ 0 \right\} \ \  \forall k>2.$$
So 
$$dim (3\mathcal{A}ss^{!})(n)=dim(3\mathcal{A}ss^{!}(2n-1))=1$$ and 
$$dim (3\mathcal{A}ss^{!})(k(n-1)+1)=0$$
for $k>2$. 
\begin{theorem}\cite{G.R}
The generating function of the dual operad of $3$-$\mathcal{A}ss$ is 
$$g_{3\mathcal{A}ss^{!}}(x)=x-x^3+x^5.$$ 
\end{theorem}
But we have seen in the previous section that any formal serie satisfying
$g_{3\mathcal{A}ss}(-s(-x))=x$ is of the form 
$$s(x)=-x+x^3-x^5+19 x^{11}+O[x]^{12}.$$
Then it doesnot correspond to the dual of $3$-$\mathcal{A}ss$.
\begin{corollary}
The quadratic operad $3\mathcal{A}ss$ is not Koszul.
\end{corollary}

\section{The quadratic operad $\widetilde{3\mathcal{A}ss}$}
In \cite{G.R.3} we have defined, given a quadratic operad $\p$, a quadratic operad $\widetilde\p$
 by the following caracteristic property:

\noindent {\it For every $\p$-algebra $A$ and every $\widetilde\p$-algebra $B$, the tensor product $A \otimes B$ 
is a $\p$-algebra.}

\begin{proposition} We have
 $$\widetilde{3\mathcal{A}ss}=3t\mathcal{A}ss.$$
\end{proposition}
Here the product is considered to be of degree $0$. In many classical cases, we have seen that $\p ^!=\widetilde\p$. Some examples where this equality 
is not realized are constructed considering binary non-associative algebras. 
\bigskip

\noindent{\bf Remark : The operad of Jordan Triple systems.}  A Jordan Triple system on a vector space is a $3$-ary product $\mu$
satisfying the commutativity condition
$$\mu(x_1,x_2,x_3)=\mu(x_3,x_2,x_1)$$
and
$$
\begin{array}{l}
 \mu(x_1,x_2,\mu(x_3,x_4,x_5))+\mu(x_3,\mu(x_2,x_1,x_4),x_5)=
\mu(\mu(x_1,x_2,x_3),x_4,x_5)\\
+\mu(x_3,x_4,\mu(x_1,x_2,x_5)).
\end{array}$$
We denote by $\mathcal{J}ord_3$ the operad of algebras defined by a Jordan triple system. In \cite{G.W} one proves that this quadratic
operad is of Koszul computing its dual. But the definition of the dual is based of the nongraded version of \cite{Gne}. Then they prove that the dual corresponds
to $3t\mathcal{A}ss$ without degree and we are not sure that the koszulity is satisfied. For instance it is not proved. 
This calculus using product of degree $1$ is making by Nicolas Goze and myself. But the delay imposed to present this paper to publication of the actes of the
Algebra, Geometry, and Mathematical Physics Workshop prevents to end this calculus.

\bigskip

\noindent{\bf Acknowledgment.} We thank the organizers of the Algebra, Geometry, and Mathematical Physics Workshop to permit
to present this work and more specially Eugen Paal and Sergei Silvestrov.


\begin{thebibliography}{9}

\bibitem{At.Mk} H. Ataguema, A. Makhlouf. Deformations of ternary algebras.
J. Gen. Lie Theory Appl. 1 (2007), no. 1, 41--55.

\bibitem{Ge} M. Gerstenhaber. On the deformation of rings and algebras. Ann.
of Math. (2) 79 1964 59--103.

\bibitem{G.K} V. Ginzburg, M. Kapranov. Koszul duality for operads.  Duke Math. J.  76  (1994),  no. 1, 203--272.



\bibitem{Gne} A. V. Gnedbaye Op\'{e}rades des alg\`{e}bres $k+1$-aires.
Operads: Proceedings of Renaissance Conferences (Hartford, CT/Luminy, 1995),
83--113, Contemp. Math., 202, Amer. Math. Soc., Providence, RI, 1997.

\bibitem{G.W} A.V. Gnedbaye,  M. Wambst. Jordan triples and operads.  J. Algebra  231  (2000),  no. 2, 744--757. 

\bibitem{G.R1} Michel Goze, Elisabeth Remm. Lie-admissible algebras and operads. J.
Algebra 273 (2004), no. 1, 129--152.

\bibitem{G.R2} Michel Goze, Elisabeth Remm. Lie admissible coalgebras.
J. Gen. Lie Theory Appl. 1 (2007), no. 1, 19--28.

\bibitem{G.R.3}  Elisabeth Remm, Michel Goze. On the algebras obtained by tensor product arXiv:math/0606105.
 
\bibitem{G.R} Nicolas Goze, Elisabeth Remm. $n$-ary  associative algebras, cohomology, 
free algebras and coalgebras. submitted to JGLTA (2008).

\bibitem{H} Eric Hoffbeck. A Poincar\'e-Birkhoff-Witt criterion for Koszul operads
arXiv:0709.2286v2 and v3. 
 


\bibitem{M.R} M. Markl, E. Remm. D\'{e}formations des op\'{e}rades non
Koszul.\ Pr\'{e}print Mulhouse 2008.

\bibitem{M.S.S} M. Markl, S. Shnider, J. Stasheff. Operads in algebra, topology and physics. 
Mathematical Surveys and Monographs, 96. American Mathematical Society, Providence, RI, 2002.


\end{thebibliography}
\end{document}